\documentclass[a4paper,11pt]{article}

\usepackage{type1cm}        
\usepackage{makeidx}         
\usepackage{graphicx}        
\usepackage{multicol}        
\usepackage[bottom]{footmisc}

\usepackage{multirow}
\usepackage{newtxtext}
\usepackage{float}
\usepackage{booktabs}

\usepackage{amsmath}
\usepackage{amssymb}

\usepackage[scale=0.8]{geometry}

\usepackage[colorlinks=true]{hyperref}
\usepackage[capitalise,nameinlink,noabbrev]{cleveref}

\makeindex             

\usepackage{nicefrac}

\begin{document}

\title{Neural network-driven domain decomposition for efficient solutions to the Helmholtz equation\thanks{\footnotesize
        V.~Dolean (ORCID 0000-0002-5885-1903, \texttt{v.dolean.maini@tue.nl})\,
        $\cdot$\,
        D.~Hrebenshchykova (ORCID 0009-0006-0190-7052, \texttt{daria.hrebenshchykova@inria.fr})\,
        $\cdot$\,
        S.~Lanteri (ORCID 0000-0001-6394-0487, \texttt{stephane.lanteri@inria.fr})\,
        $\cdot$\,
        V.~Michel-Dansac (ORCID 0000-0002-3859-8517, \texttt{victor.michel-dansac@inria.fr})}}

\author{
    Victorita Dolean\textsuperscript{1} \quad
    Daria Hrebenshchykova\textsuperscript{2} \quad
    Stéphane Lanteri\textsuperscript{2} \quad
    Victor Michel-Dansac\textsuperscript{3} \\[8pt]
    \small
    \textsuperscript{1}Department of Mathematics and Computer Science, Eindhoven University of Technology, The Netherlands \\
    \small
    \textsuperscript{2}Atlantis project-team, Inria Research Center at Université Côte d'Azur, France \\
    \small
    \textsuperscript{3}MACARON project-team, Université de Strasbourg, CNRS, Inria, IRMA, France
}

\date{February 04, 2026}

\maketitle

\begin{abstract}
    Accurately simulating wave propagation is crucial in fields such as acoustics, electromagnetism, and seismic analysis. Traditional numerical methods, like finite
    difference and finite element approaches, are widely used to solve governing partial differential equations (PDEs) such as the Helmholtz equation. However, these methods face significant computational challenges when applied to high-frequency wave problems in complex two-dimensional domains.
    This work investigates Finite Basis Physics-Informed Neural Networks (FBPINNs) and their multilevel extensions as a promising alternative. These methods leverage domain decomposition, partitioning the computational domain into overlapping subdomains, each governed by a local neural network. We assess their accuracy and computational efficiency in solving the Helmholtz equation for the homogeneous case, demonstrating their potential to mitigate the limitations of traditional approaches.
\end{abstract}

\section{Introduction}
\label{sec:1}
Accurate simulation of wave propagation is crucial in many scientific and engineering applications, where the frequency-domain Helmholtz equation plays a central role. A particular challenge arises when solving this equation in free space, since the correct mathematical formulation requires the Sommerfeld radiation condition at infinity. This condition can be approximated by an absorbing boundary, such as a Robin boundary condition or perfectly matched layers (PMLs). Handling such open-domain problems accurately is essential to avoid reflections and ensure physically consistent wavefields.

We focus on solving a well-posed time-independent problem, in the form
\begin{equation}
    \mathcal{D}[u](\mathbf{x}) = g(\mathbf{x}), \quad \mathbf{x} \in \Omega \subset \mathbb{R}^d,
    \label{eq:main_system}
\end{equation}
where \( \mathcal{D} \) is some differential operator, \( u \) is the unknown solution and  $\mathbf{x}$ is the space variable, with $\Omega$ being the space domain. Note that \eqref{eq:main_system} does not contain a boundary operator; this is intentional as boundary conditions will be handled by PMLs.
To solve~\eqref{eq:main_system}, we use Physics-Informed Neural Networks (PINNs) proposed in~\cite{pinns}. PINNs train a Neural Network (NN) to solve a differential equation  by directly approximating the solution, i.e., defining \( v(\cdot, \boldsymbol{\theta}) \approx u(\cdot) \), where $\boldsymbol{\theta}$ represents a set of $P$ trainable parameters of the neural network.
The main advantages of PINNs are their mesh-free nature and the fact that we end up with
a continuous functional approximation of the solution, defined over the entire domain, rather than a discretised solution on a fixed mesh as in conventional numerical methods.
Based on~\cite{pinns}, the following loss function is minimized to train
the PINN
\begin{equation}
    \begin{aligned}
        L(\boldsymbol{\theta}) & = \frac{1}{N_I} \sum_{i=1}^{N_I}
        \underbrace{
            \big( \mathcal{D}[v_{\boldsymbol{\theta}}](\mathbf{x}_i) - g(\mathbf{x}_i) \big)^2
        }_{\text{PDE residual}}
        ,
    \end{aligned}
    \label{eq:pinn_loss}
\end{equation}
where \smash{\( \{\mathbf{x}_i\}_{i=1}^{N_I} \)} is a set of collocation points sampled in the domain, and where we have defined $v_{\boldsymbol{\theta}} = \mathbf{x} \mapsto v(\mathbf{x}, \boldsymbol{\theta})$. As mentioned above, we do not include a boundary loss term here, as the PML formulation ensures wave absorption at the boundaries, removing the necessity for explicit boundary constraints as in standard PINNs.

We move from the standard PINNs formulation to its domain-decomposed extension, namely the Finite Basis PINNs (FBPINNs) introduced in~\cite{fbpinns}, which decomposes the domain into overlapping subdomains. In this paper, we combine FBPINNs with PMLs to simulate free-space wave propagation. The goal of this work is thus to assess the performance of the combination of FBPINNs and PMLs for solving the Helmholtz equation.

\section{Mathematical Model}
\label{sec:2}
In this work, we consider the frequency-domain Helmholtz equation in free space, given by
\begin{equation*}
    -\Delta u(\mathbf{x}) - k^2 u(\mathbf{x}) = g(\mathbf{x}),\quad \mathbf{x} \in \Omega,
    \label{eq:helmholtz}
\end{equation*}
for $\mathbf{x} = (x, y)^\intercal \in \Omega = \mathbb{R}^d$ with $d=2$. Here \smash{$ k = \nicefrac{2\pi}{\lambda}$} is the constant (scalar) wavenumber defined by wavelength \smash{$\lambda = \nicefrac c f$}, where $c$ is the wave speed and $f$ is the frequency.
The unknown function is $u: \mathbb{R}^2 \to \mathbb{C}$ and $g: \mathbb{R}^2 \to \mathbb{C}$ is a source function.
This problem naturally falls into the general framework \eqref{eq:main_system}.
Physically valid solutions satisfy the Sommerfeld radiation condition at infinity
\begin{equation*}
    \lim_{|\mathbf{x}| \to \infty} |\mathbf{x}|^{\frac{d-1}{2}} \left( \frac{\partial u}{\partial |\mathbf{x}|} - i k u \right) = 0.
\end{equation*}

\subsection{Perfectly Matched Layer}

To approximate the Sommerfeld radiation condition, the problem is reformulated using a PML transformation introduced by Bérenger \cite{pml_basic} in 1994. The PML offers a more robust approximation than e.g.\ Robin boundary conditions by creating a finite, artificial absorbing layer around the simulation area.

\begin{figure}[!ht]
    \centering
    \includegraphics[width=0.6\linewidth]{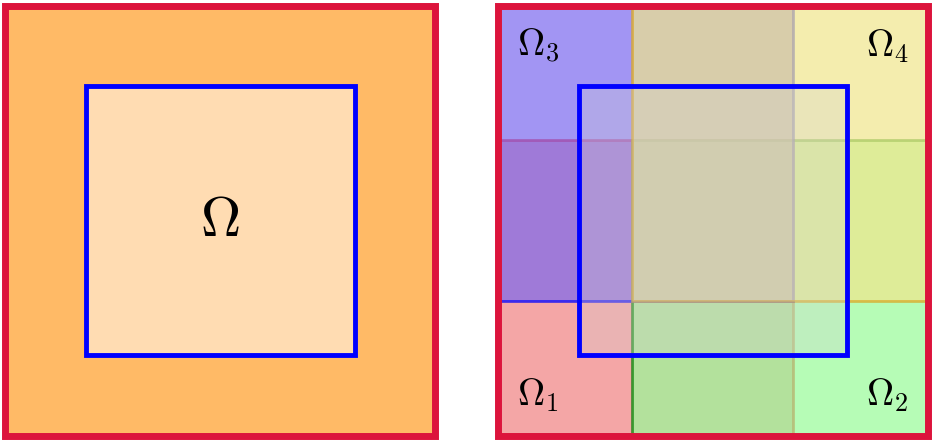}
    \caption{
        Schematic illustration of the computational domains for PINNs (left) and FBPINNs (right) with PML.
        The blue rectangle denotes the physical domain $\Omega$ where we want to find the solution.
        The red line indicates the outer boundary of the computational domain, including the PML region.
        The PML is thus located between the red and blue boundaries and ensures wave absorption.
        In the FBPINNs case (right), the physical domain together with the PML region is further divided into four overlapping subdomains  $\smash{\{\Omega_j\}_{j \in \{1, \dots, 4\}}}$.
    }
    \label{fig:pml_domains}
\end{figure}
We are looking for the solution in the region $\Omega = [-L,L]^2$, considering a symmetrical rectangular domain without loss of generality, since equivalent results can be obtained for non-symmetric rectangular domains. While the method can be adapted to more general domain geometries by modifying the functions defined below, such extensions are not considered in this work for the sake of conciseness. To apply the PML, the computational domain  $\Omega= [-L-L_\text{PML},L + L_\text{PML}]^2$, where $L_\text{PML}$ represents the width of our absorbing layer.
The computational domain is illustrated in \cref{fig:pml_domains} (left), where the blue square corresponds to the physical region and the red layer represents the PML.  At this stage, we consider a single computational domain. The implementation of a PML uses a coordinate stretching technique, thus, we define the coordinate stretching functions
\begin{equation*}
    D_x: x \mapsto \frac{1}{1 - i \frac{\sigma(x)}{k}}
    \text{\quad and \quad}
    D_y: y \mapsto \frac{1}{1 - i \frac{\sigma(y)}{k}}.
\end{equation*}
where
\[
    \sigma(x)=\sigma_{0}k
    \left(\dfrac{\max\!\left(|x|-L,\,0\right)}{L_{\text{PML}}}\right)^2
\]

The tuning constant $\sigma_{0}$ and the PML thickness are specified parameters in the numerical experiments. The function $\sigma$ is defined following \cite{pml_dd, pml_ref}, but other choices are also possible.
Then, thanks to the absorbing layer, we no longer impose a boundary condition explicitly, since it is incorporated into the equation itself. As a result, the differential operator with PML $\mathcal{D}_\text{PML}$, associated with the Helmholtz equation, is defined as
\begin{equation}
    \label{eq:PML_eq}
    \mathcal{D}_\text{PML} = -D_x\frac{\partial}{\partial x} \left( D_x \frac{\partial  }{\partial x} \right) - D_y\frac{\partial}{\partial y} \left( D_y \frac{\partial  }{\partial y}\right) - k^2.
\end{equation}

\subsection{FBPINNs}
This section defines the mathematical formulation of FBPINNs according to \cite{fbpinns}.
First, we split the global domain $\Omega$ into $J$ overlapping subdomains
$\smash{\{\Omega_j\}_{j \in \{1, \dots, J\}}}$. To achieve smoother and more accurate solutions, these subdomains are constructed with an overlap $\eta$, which will be specified in the numerical experiments. This overlap is introduced to avoid dealing explicitly with continuity constraints at subdomain interfaces.
The approximate solution to our problem \eqref{eq:main_system} is then defined as
\begin{equation}
    \tilde{u}(\mathbf{x}, \boldsymbol{\theta})
    = \sum_{j=1}^{J} \varphi_j(\mathbf{x}) \, v_j(\mathbf{x}; \boldsymbol{\theta}_j),
    \label{eq:fbpinn_aprox}
\end{equation}
where $\tilde{u}$ is the weighted combination of the local neural networks $v_j$,
each trained on its corresponding subdomain $\Omega_j$.
The window function $\varphi_j$ restricts the influence of~$v_j$
to its local subdomain while ensuring a smooth transition between neighboring regions.

Each window function must be smooth, differentiable, and form a partition of unity, i.e.
$\smash{\sum_{j=1}^{J} \varphi_j(\mathbf{x}) = 1}$.
Like in \cite{multilevel}, we use a cosine-based window function
\begin{equation*}
    \varphi_j(\mathbf{x})
    = \frac{\hat{\varphi}_j(\mathbf{x})}
    {\sum_{m=1}^{J} \hat{\varphi}_m(\mathbf{x})},
    \quad \text{where} \quad
    \hat{\varphi}_j(\mathbf{x})
    = \prod_{n=1}^{2}\left[1 + \cos\!\left(\pi \frac{x_n - \mu_{j,n}}{\psi_{j,n}}\right)\right]^2,
    \label{eq:window_function}
\end{equation*}
truncated to $0$ outside the subdomain $\Omega_j$.
Here, $\mu_{j,n}$ denotes the center of the $j$-th subdomain along the $n$-th spatial axis,
and $\smash{\psi_{j,n} = \nicefrac{\eta \Delta_n }{2}}$ defines its half-width, with
$\Delta_n$ the subdomain size along dimension $n$.

The FBPINNs loss function is based on \eqref{eq:pinn_loss}, but now combined with PML, we no longer include the boundary residual term, since the PML formulation ensures wave absorption at the boundaries. Consequently, the computational domain becomes the extended region $\Omega_{\text{PML}}$, and the standard differential operator $\mathcal{D}$ is replaced by the PML operator $\mathcal{D}_\text{PML}$ defined in \eqref{eq:PML_eq}.
The general structure of the FBPINNs combined with PML is shown in \cref{fig:pml_domains} (right), where the global domain, including the PML, is partitioned into four overlapping subdomains.
Given the ansatz defined by \eqref{eq:fbpinn_aprox}, the loss function accounting for both FBPINNs and PML takes the following form
\begin{equation*}
    L(\boldsymbol{\theta}) = \frac{1}{N} \sum_{i=1}^{N} \big( \mathcal{D}_\text{PML} [ \tilde{u}_{\boldsymbol{\theta}}](\mathbf{x}_i) - g(\mathbf{x}_i) \big)^2,
    \label{eq:fbpinn_loss}
\end{equation*}
where $\tilde{u}_{\boldsymbol{\theta}} = \mathbf{x} \mapsto \tilde{u}(\mathbf{x}, \boldsymbol{\theta})$.

\section{Numerical Results}
\label{sec:3}

In this section we present numerical results for the 2D Helmholtz problem with a Gaussian point source $g(x,y) = e^{-k(x^2 + y^2)}$ in the domain $\Omega = [-5,5]^2$.
In our tests, we set the wave speed $c = 1$ and vary the frequency $f$ such that the wavenumber $\smash{k=\nicefrac{2\pi}{\lambda}}$ does not coincide with the eigenfrequencies of the domain, ensuring the existence and uniqueness of the solution.
We vary the PML width based on the wavelength $\smash{\lambda = \nicefrac{2\pi}{k}}$, specifically $0.25\lambda, 0.5\lambda \text{ and } \lambda$. To ensure that subdomains overlap, the overlap ratio must satisfy $\eta>1$. If it is smaller than 1, the subdomains remain disjoint. In our tests we use an overlap $\eta=1.5$.
To help the network find the correct solution for the imaginary part of the problem, we start training with an incorrect source term for the imaginary part that contains small oscillations with wavenumber $k$. After a few epochs, we replace the source term with the correct one, which is zero in our case.

Some PINNs and FBPINNs settings and hyperparameters are fixed across all experiments. Specifically, we use a learning rate of 0.001, a standard multilayer perceptron, 5000 collocation points and a \textit{sine} activation function. We divide the domain into four subdomains when using FBPINNs.

We compare the Adam + L-BFGS optimization technique with the Energy Natural Gradient Descent (ENGD) method. Adam and L-BFGS are state-of-the-art optimizers for training PINNs, commonly used in the literature due to their reliability and convergence properties.
ENGD, introduced in \cite{NG}, arises by projecting the optimal function-space descent direction onto the tangent space of the neural-network manifold, yielding a preconditioned gradient in parameter space, that best approximates this optimal direction. It modifies the standard gradient update as follows:
\begin{equation*}
    \boldsymbol{\theta}_{t+1} = \boldsymbol{\theta}_t - \, G^\dagger(\boldsymbol{\theta}_t)\, \nabla_\theta L(\boldsymbol{\theta}_t),
\end{equation*}
where $G^\dagger$ is the Moore-Penrose pseudoinverse of the matrix $G$ with components
\[
    \big(G(\boldsymbol{\theta})\big)_{ij} = \int_{\Omega_\text{PML}} \,
    \partial_{\theta_i} \! \big(\mathcal{D}_\text{PML}[\tilde u_{\boldsymbol{\theta}}]\big) \! (\mathbf{x}) \
    \partial_{\theta_j} \! \big(\mathcal{D}_\text{PML}[\tilde u_{\boldsymbol{\theta}}]\big) \! (\mathbf{x}) \
    d\mathbf{x}
\]
for all $(i,j) \in \{1, \dots, P\}$.
For experiments with ENGD and Adam optimizer in combination with L-BFGS, where L-BFGS is initiated after $70\%$ of the total training epochs, we use the architecture 3 hidden layers with 32 neurons each for PINNs and 3 hidden layers with 16 neurons in every layer for each NN of FBPINNs.

To evaluate the effect of domain decomposition, we compare the performance of standard PINNs and FBPINNs for different wavenumbers and PML thicknesses.
In \cref{tab:l2_real,tab:l2_imag}, we report the relative $L^2$ error, defined by
\begin{equation*}
    \varepsilon_{L^2} =
    \frac{\| u_{\mathrm{approx}} - u_{\mathrm{ref}} \|_2}
    {\| u_{\mathrm{ref}} \|_2},
    \label{eq:relL2}
\end{equation*}
where $u_{\mathrm{approx}}$ and $u_{\mathrm{ref}}$ represent the approximate and reference solutions, respectively. We compute the errors for the real and the imaginary parts separately.
\begin{table}[!ht]
    \centering
    \setlength{\tabcolsep}{2.5pt}
    \renewcommand{\arraystretch}{0.8}

    \begin{tabular}{cc}
        \begin{minipage}{0.48\textwidth}
            \centering
            \textbf{(a) Adam + L-BFGS optimizer}\\[2pt]
            \begin{tabular}{lccccc}
                \toprule
                \textbf{Method} & \textbf{$k$} & \textbf{0.25$\lambda$} & \textbf{0.5$\lambda$} & \textbf{1.0$\lambda$} \\
                \midrule
                PINNs           & 0.57         & 5.6e--2                & 1.6e--2               & 4.0e--3               \\
                                & 1.59         & 4.0e--1                & 2.9e--2               & 2.2e--2               \\
                                & 4.51         & 9.1e--1                & 9.7e--1               & 9.7e--1               \\
                \cmidrule(lr){1-5}
                FBPINNs         & 0.57         & 4.1e--2                & 1.5e--2               & 7.3e--3               \\
                                & 1.59         & 3.0e--1                & 1.6e--2               & 1.1e--2               \\
                                & 4.51         & 9.3e--1                & 9.3e--1               & 9.5e--1               \\
                \bottomrule
            \end{tabular}
        \end{minipage}
         &
        \begin{minipage}{0.48\textwidth}
            \centering
            \textbf{(b) ENGD optimizer}\\[2pt]
            \begin{tabular}{lccccc}
                \toprule
                \textbf{Method} & \textbf{$k$} & \textbf{0.25$\lambda$} & \textbf{0.5$\lambda$} & \textbf{1.0$\lambda$} \\
                \midrule
                PINNs           & 0.57         & 1.7e--1                & 1.0e--2               & 5.4e--3               \\
                                & 1.59         & 6.7e--1                & 1.8e--2               & 4.9e--3               \\
                                & 4.51         & 6.7e--1                & 3.0e--1               & 1.4e--2               \\
                \cmidrule(lr){1-5}
                FBPINNs         & 0.57         & 7.3e--1                & 6.1e--3               & 3.3e--3               \\
                                & 1.59         & 8.7e--1                & 1.5e--2               & 4.1e--3               \\
                                & 4.51         & 6.4e--1                & 1.7e--1               & 1.1e--2               \\
                \bottomrule
            \end{tabular}
        \end{minipage}
        \\
    \end{tabular}
    \caption{Relative $L_2$ error between real part of approximate and reference solutions for different $L_{\mathrm{PML}}$, comparing PINNs and FBPINNs performance using (a) Adam + L-BFGS and (b) ENGD.}
    \label{tab:l2_real}
\end{table}


\begin{table}[!ht]
    \centering
    \setlength{\tabcolsep}{2.5pt}
    \renewcommand{\arraystretch}{0.8}

    \begin{tabular}{cc}
        \begin{minipage}{0.48\textwidth}
            \centering
            \textbf{(a) Adam + L-BFGS optimizer}\\[2pt]
            \begin{tabular}{lccccc}
                \toprule
                \textbf{Method} & \textbf{$k$} & \textbf{0.25$\lambda$} & \textbf{0.5$\lambda$} & \textbf{1.0$\lambda$} \\
                \midrule
                PINNs           & 0.57         & 1.0e--1                & 2.0e--2               & 6.2e--3               \\
                                & 1.59         & 6.4e--1                & 2.7e--1               & 8.8e--1               \\
                                & 4.51         & 1.0e+0                 & 1.0e+0                & 1.0e+0                \\
                \cmidrule(lr){1-5}
                FBPINNs         & 0.57         & 5.6e--1                & 2.7e--2               & 1.8e--3               \\
                                & 1.59         & 5.1e--1                & 5.9e--2               & 3.3e--2               \\
                                & 4.51         & 1.0e+0                 & 1.0e+0                & 1.0e+0                \\

                \bottomrule
            \end{tabular}
        \end{minipage}
         &
        \begin{minipage}{0.48\textwidth}
            \centering
            \textbf{(b) ENGD optimizer}\\[2pt]
            \begin{tabular}{lccccc}
                \toprule
                \textbf{Method} & \textbf{$k$} & \textbf{0.25$\lambda$} & \textbf{0.5$\lambda$} & \textbf{1.0$\lambda$} \\
                \midrule
                PINNs           & 0.57         & 2.7e--1                & 1.4e--2               & 1.8e--3               \\
                                & 1.59         & 8.5e--1                & 2.4e--2               & 1.4e--3               \\
                                & 4.51         & 9.0e--1                & 4.7e--1               & 8.2e--2               \\
                \cmidrule(lr){1-5}
                FBPINNs         & 0.57         & 5.9e--1                & 1.3e--2               & 1.3e--3               \\
                                & 1.59         & 6.1e--1                & 2.0e--2               & 1.3e--3               \\
                                & 4.51         & 9.0e--1                & 3.7e--1               & 1.6e--2               \\
                \bottomrule
            \end{tabular}
        \end{minipage}
        \\
    \end{tabular}
    \caption{Relative $L_2$ error between imaginary part of approximate and reference solutions for different $L_{\mathrm{PML}}$, comparing PINNs and FBPINNs performance using (a) Adam + L-BFGS and (b) ENGD.}
    \label{tab:l2_imag}
\end{table}

\begin{figure}[!ht]
    \centering

    \includegraphics[width=0.7\textwidth]{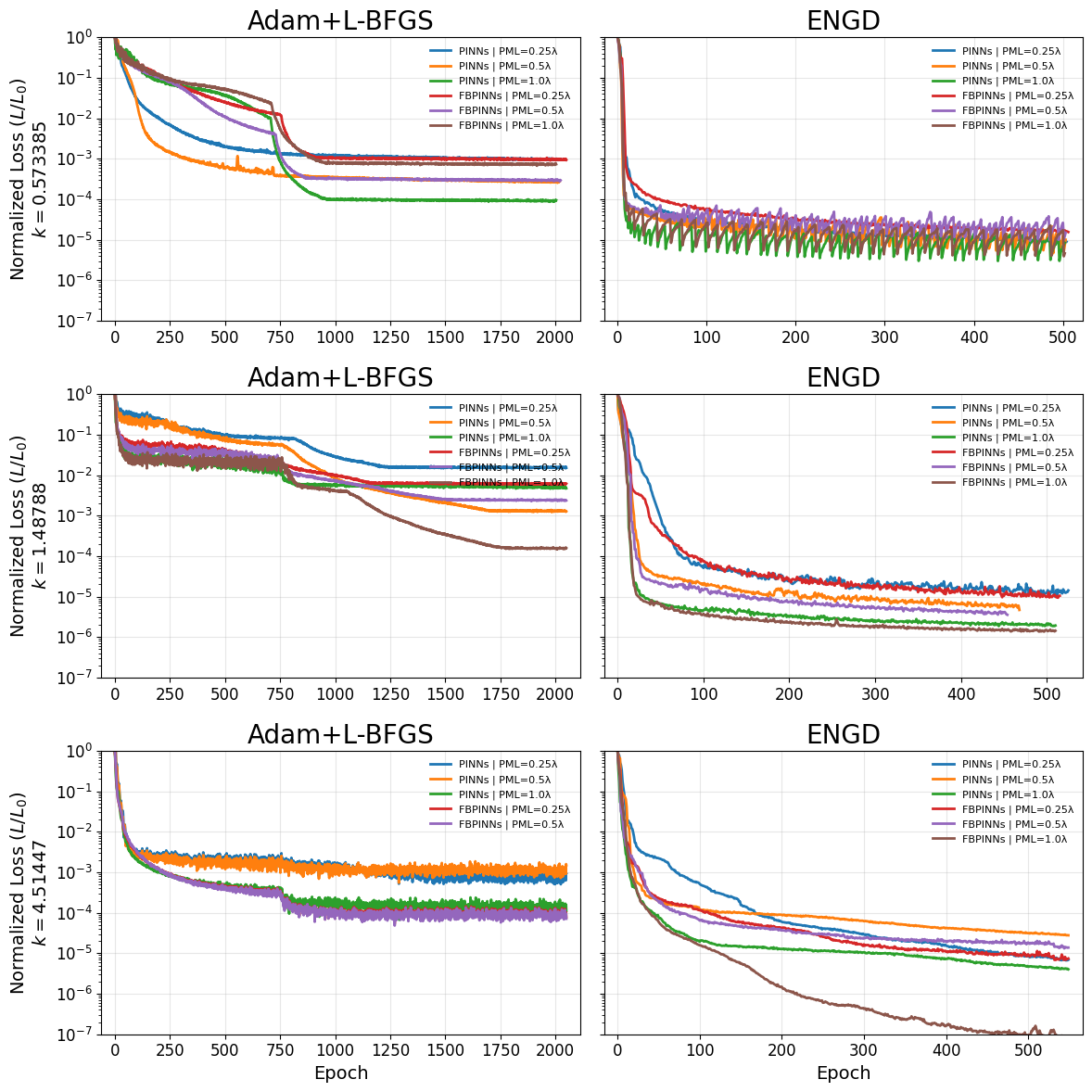}
    \caption{Normalized physics loss during training for Adam + L-BFGS (left) and ENGD (right).
        Rows correspond to different $k$ and curves compare PINNs vs FBPINNs across PML thicknesses.
    }
    \label{fig:loss_all_methods}
\end{figure}

From \cref{tab:l2_real} we can see that both optimization methods perform similarly at low frequencies ($k=0.57$). However, as the wavenumber increases, the standard PINNs quickly deteriorate, showing errors close to unity already for $k=4.51$, while the FBPINNs maintain a slightly better accuracy.  We can observe a similar situation with the imaginary part from \cref{tab:l2_imag}, but by comparing the results, we can conclude that despite our assistance to the imaginary part in the form of preliminary training with another source term, it is still approximated worse than the real part.  We assume that this effect appears because there is a local minimum where the imaginary part tends to vanish, which makes the optimization harder for this part.
The difference between PINNs and FBPINNs is clearly visible in the convergence curves in \cref{fig:loss_all_methods}, where the FBPINNs show a faster decrease in loss and achieve a lower final value, especially for high frequencies.
Nevertheless, when the width of PML is too thin (for example $L_\text{PML} = 0.25\lambda$), both approaches struggle to converge due to insufficient wave absorption at the boundaries.

\begin{figure}[!ht]
    \centering
    \includegraphics[width=0.9\textwidth]{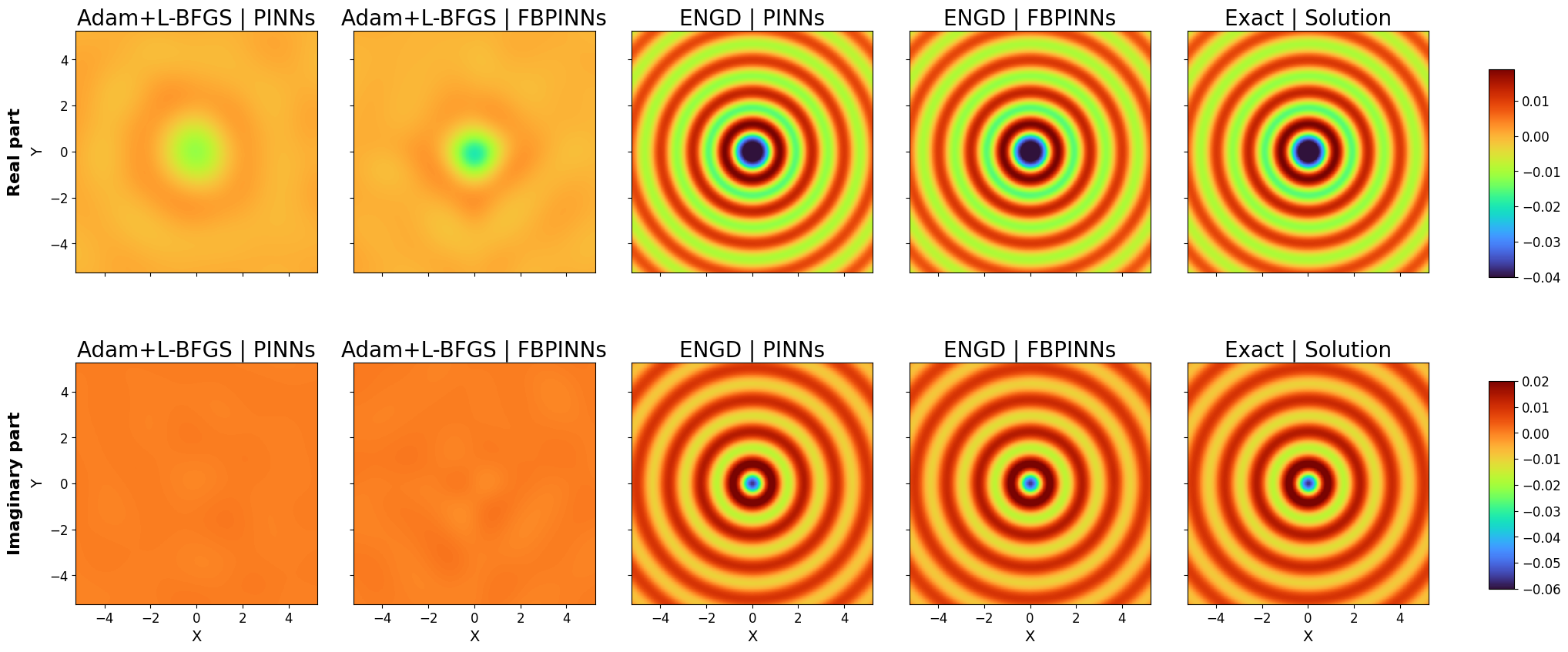}
    \caption{Approximation and reference solutions of the real and imaginary parts for the Helmholtz problem with $k=4.51$ comparing Adam+L-BFGS and ENGD optimization methods for PINNs and FBPINNs.}
    \label{fig:results}
\end{figure}

We also compare the performance of Adam + L-BFGS to ENGD.
At low frequencies, both optimizers provide similar results, although ENGD  demonstrates faster stabilization of the loss that we can see from \cref{fig:loss_all_methods}.
As the wavenumber increases, the difference becomes clearer: while the standard optimizer struggles to converge,  ENGD succeeds in maintaining stable training and correctly reproducing the oscillatory solution.
From the real and imaginary part approximations shown in \cref{fig:results}, we observe that ENGD captures the fine wave structures missed by Adam + L-BFGS, even in the imaginary part.
Regarding the errors reported in \cref{tab:l2_real,tab:l2_imag}, the results are most evident when we move to the wavenumber $k = 4.5$, since in this case, especially for the imaginary part, Adam + L-BFGS gives the worst error, while ENGD achieves an accurate approximation.
The training speed of ENGD is about six times slower than that of Adam + L-BFGS for both PINNs and FBPINNs. This is mainly due to the additional computation of the energy matrix at each iteration, which significantly increases the computational cost.
The use of double precision arithmetic further contributes to the slowdown. In practice, switching to hardware that handles double precision more efficiently could reduce this gap, although the method would still remain slower overall.

\section{Conclusion}
In this work, we showed that FBPINNs combined with PML boundaries can effectively handle wave absorption and provide more accurate solutions than standard PINNs. While the improvement remains moderate, the obtained results confirm that the approach is effective and has potential for further refinement. The use of natural gradient preconditioning also helped to stabilize training and achieve better convergence. In future work, we plan to refine the combination of FBPINNs and PML by introducing internal absorbing layers and to extend the framework toward multilevel FBPINNs as in \cite{multilevel} for more complex wave problems (higher frequencies, heterogeneous media, \ldots).

\section*{Acknowledgements}
As part of the ``France 2030'' initiative, this work has benefited from a national grant managed by the French National Research Agency (Agence Nationale de la Recherche) attributed to the Exa-MA project of the NumPEx PEPR program, under the reference ANR-22-EXNU-0002.

\section*{Competing Interests}
The authors have no conflicts of interest to declare that are relevant to the content of this chapter.

%
%
%


\bibliographystyle{plain}
\bibliography{./bibliography}

@article{pml_basic,
title = {A perfectly matched layer for the absorption of electromagnetic waves},
journal = {Journal of Computational Physics},
volume = {114},
number = {2},
pages = {185-200},
year = {1994},
issn = {0021-9991},
doi = {https://doi.org/10.1006/jcph.1994.1159},
url = {https://www.sciencedirect.com/science/article/pii/S0021999184711594},
author = {Jean-Pierre Berenger},
}

@inbook{pml_dd,
   title={Numerical Assessment of PML Transmission Conditions in a Domain Decomposition Method for the Helmholtz Equation},
   ISBN={9783031507694},
   ISSN={2197-7100},
   url={http://dx.doi.org/10.1007/978-3-031-50769-4_53},
   DOI={10.1007/978-3-031-50769-4_53},
   booktitle={Domain Decomposition Methods in Science and Engineering XXVII},
   publisher={Springer Nature Switzerland},
   author={Bootland, Niall and Borzooei, Sahar and Dolean, Victorita and Tournier, Pierre-Henri},
   year={2024},
   pages={445–453} }

@article{multilevel,
   title={Multilevel domain decomposition-based architectures for physics-informed neural networks},
   volume={429},
   ISSN={0045-7825},
   url={http://dx.doi.org/10.1016/j.cma.2024.117116},
   DOI={10.1016/j.cma.2024.117116},
   journal={Computer Methods in Applied Mechanics and Engineering},
   publisher={Elsevier BV},
   author={Dolean, Victorita and Heinlein, Alexander and Mishra, Siddhartha and Moseley, Ben},
   year={2024},
   month=sep, pages={117116} }

@article{fbpinns,
   title={Finite basis physics-informed neural networks ({FBPINNs}): a scalable domain decomposition approach for solving differential equations},
   volume={49},
   ISSN={1572-9044},
   url={http://dx.doi.org/10.1007/s10444-023-10065-9},
   DOI={10.1007/s10444-023-10065-9},
   number={4},
   journal={Advances in Computational Mathematics},
   publisher={Springer Science and Business Media LLC},
   author={Moseley, Ben and Markham, Andrew and Nissen-Meyer, Tarje},
   year={2023},
   month=jul }

@inproceedings{NG,
  author    = {M\"{u}ller, Johannes and Zeinhofer, Marius},
  booktitle = {Proceedings of the 40th International Conference on Machine Learning},
  title     = {{Achieving High Accuracy with PINNs via Energy Natural Gradient Descent}},
  year      = {2023},
  editor    = {Krause, A. and Brunskill, E. and Cho, K. and Engelhardt, B. and Sabato, S. and Scarlett, J.},
  pages     = {25471--25485},
  publisher = {PMLR},
  series    = {Proceedings of Machine Learning Research},
  volume    = {202},
  _month_   = {23--29 Jul},
  pdf       = {https://proceedings.mlr.press/v202/muller23b/muller23b.pdf},
  url       = {https://proceedings.mlr.press/v202/muller23b.html}
}

@article{pinns,
title = {Physics-informed neural networks: A deep learning framework for solving forward and inverse problems involving nonlinear partial differential equations},
journal = {Journal of Computational Physics},
volume = {378},
pages = {686-707},
year = {2019},
issn = {0021-9991},
doi = {https://doi.org/10.1016/j.jcp.2018.10.045},
url = {https://www.sciencedirect.com/science/article/pii/S0021999118307125},
author = {Maziar Raissi and Paris Perdikaris and George Em Karniadakis},
}

@article{pml_ref,
author = {Song, Chao and Alkhalifah, Tariq and Bin Waheed, Umair},
year = {2021},
month = {10},
title = {A versatile framework to solve the {Helmholtz} equation using physics-informed neural networks},
volume = {228},
journal = {Geophysical Journal International},
doi = {10.1093/gji/ggab434}
}

\end{document}